\documentclass{article}
\usepackage{arxiv}

\usepackage[utf8]{inputenc}
\usepackage[T1]{fontenc}
\usepackage[ukrainian,russian,english]{babel}
\usepackage[unicode, pdftex]{hyperref}
\usepackage{url}            
\usepackage{booktabs}       
\usepackage{amsfonts}       
\usepackage{nicefrac}       
\usepackage{microtype}      
\usepackage{lipsum}		
\usepackage{graphicx}
\usepackage{amsmath}
\usepackage{amsfonts}
\usepackage{amssymb}  
\usepackage{natbib}
\usepackage{doi}

\title{Algorithms and topological invariants for dynamic systems. I. Basic definitions}

\author{ \href{https://orcid.org/0000-0002-7164-807X}{\includegraphics[scale=0.06]{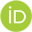}
\hspace{1mm}Alexandr O.~Prishlyak}\thanks{https://sites.google.com/view/prof-prishlyak, https://orcid.org/0000-0002-7164-807X} \\
	Department of Computer Methods of \\ 
	Mechanics and Control Processes\\
	Taras Shevchenko National University of Kyiv\\
	Kyiv, Ukraine \\
	\texttt{prishlyak@knu.ua} }

\renewcommand{\shorttitle}{Algorithms and topological invariants for dynamic systems. I. Basic definitions}

\newtheorem{theorem}{Theorem}

\newtheorem{lema}{Lemma}

\hypersetup{
pdftitle={Algorithms and topological invariants for dynamic systems. I. Basic definitions},
pdfauthor={Alexandr O. Prishlyak},
pdfkeywords={Topological classification, Morse function, Morse-Smale flow},
}

\begin{document}
\maketitle

\begin{abstract}
We construct algorithms and topological invariants that allow us to distinguish the topological type of a surface, as well as functions and vector fields for their topological equivalence. 
In the first part we discus the main structures used in the topology of manifolds: vector fields, dynamical systems, Morse functions, cell decompositions, and the fundamental group.
\end{abstract}

\keywords{Topological classification \and Morse function \and Morse-Smale flow}

\section*{Introduction}
Topological properties of manifolds, functions, and dynamical systems are often studied through the construction of topological invariants that have a discrete nature, meaning they can be described using a finite set of integers, which allows for the use of computational techniques when working with them. Morse theory enables the construction of a cell complex structure on a manifold; however, continuous mappings used to attach cells are generally difficult to encode. Therefore, triangulations are more commonly used in computer modeling. Their drawback is the large number of simplices required to construct such structures. A compromise solution between these two structures is regular cell complexes. To find all possible structures under study, efficient algorithms for their recognition are necessary. Topological invariants will be useful if there are efficient algorithms for their computation and comparison. The issues mentioned are addressed in low dimensions (up to four) in this tutorial.

The first part discusses the main structures used in the topology of manifolds: vector fields, dynamical systems, Morse functions, cell decompositions, and the fundamental group.

The second part examines discrete structures for which computers can be used for their specification and manipulation.

The third part describes algorithms that allow for the recognition of manifolds and some of their properties for the structures described in the second chapter.

The fourth part is dedicated to describing possible topological structures of functions and dynamical systems on low-dimensional manifolds.


\newpage

\section{Basic Topological Definitions and Notations}

In this section, we will recall the concept of topological structure. It is the most general structure that allows defining a continuous mapping.

\textbf{Definition.}   \index{topological structure}
\textit{A topological structure} (or \textit{topology}) on a set $X$ is defined as a collection of subsets called open sets, which have the following properties:

1) The empty set and the set $X$ are open,

2) The union of any number of open sets is an open set,

3) The intersection of a finite number of open sets is an open set.

Elements of $X$ are called points, and the set on which the topology is defined is called a topological space.

We consider subsets of Euclidean space, which is the n-dimensional vector space $\mathbb{R}^n$ with the standard (Euclidean) distance defined as $$d(x,y)=\sqrt{\sum_{i=1}^n (x_i-y_i)^2 }$$ for points $ x=(x_1,x_2,\ldots ,x_n), y=(y_1,y_2,\ldots , y_n).$

An $a$-neighborhood of a point $x$ ($a>0$) is defined as the set
$$ B(x,a)= \{ y : d(x,y) <a\}.$$

For subsets of Euclidean (as well as metric) space, the standard topology is defined as follows: a subset $U$ is open if for every point $x\in U$, there exists an $a>0$ such that $B(x,a) \subset U$.

\textbf{Definition.}   \index{closed set}
A subset $A \subset X$ is called \textit{closed} if its complement $X \ A$ is an open set.

\textbf{Definition.}   \index{neighborhood}
An \textit{open neighborhood} of a point is any open set that contains it. A \textit{neighborhood} (in a broader sense) is any set that contains an open neighborhood.

\textbf{Definition.}   \index{continuous mapping}
A \textit{continuous mapping} $f:X\to Y$ between topological spaces is a mapping for which the preimage of any open set from $Y$ is open in $X$. The set of all continuous mappings from $X$ to $Y$ is denoted by $C(X,Y)$. \index{homeomorphism} A \textit{homeomorphism} is a continuous bijective mapping such that the inverse mapping is also continuous.

\subsection{Topological Properties}

\textit{Topological properties} are properties that are preserved under homeomorphisms.

\textbf{Definition.}   \index{connected space}
A space is called \textit{connected} if it cannot be represented as the union of two non-empty, disjoint open sets. A space $X$ is \textit{linearly connected} if for any points $x,y \in X$ there exists a path $\varphi$ connecting them ($\varphi \in C([0,1], X), \varphi (0)=x, \varphi(1)=y$.) From linear connectivity follows the property of \textit{connected component}, which is the maximal connected set with respect to inclusion.

\textbf{Definition.}   \index{Hausdorff space}
A space is called \textit{Hausdorff} if for any two points in it there exist neighborhoods that do not intersect. An open cover of a space is a collection of open sets whose union gives the entire space. \index{compact space} A space is \textit{compact} if from any open cover of it, a finite subcollection can be extracted that forms an open cover of the space. Subsets in Euclidean space are compact if they are closed and bounded (there exists a number such that the distance between any two points is limited by this number).

They say that a property holds locally in the space $X$ if for every point $x \in X$, there exists a neighborhood where this property holds. In particular, the space is called \textit{locally connected} if there exists a connected neighborhood for every point. The space is \textit{locally compact} if for every point there is a neighborhood contained within a compact set. The space is \textit{locally Euclidean} if for every point there is a neighborhood homeomorphic to Euclidean space.

\subsection{Topological Operations}

\textbf{Definition.}   \index{subspace}
A \textit{subspace} $A$ of the space $X$ is a subset $A \subset X$ that has a topological structure defined on it, where the elements are the intersections of $A$ with open sets in $X$.

\textbf{Definition.}   \index{factor space}
Let there be an equivalence relation $\sim$ on the space $X$, $[x]$ is the equivalence class of the point $x \in X$, and $Y = X / \sim$ is the set of equivalence classes, with $p: X \to Y$ being the projection, where $p(x) = [x]$. Then, in $Y = X / \sim$, the \textit{factor topology} is defined such that a set $U$ is open if $p^{-1}(U)$ is open in $X$. The factor set $X / \sim$ together with the factor topology is called a \textit{factor space}.

The factor space of the space $X$ with respect to the subspace $A \subset X$ is called the factor space with respect to the equivalence relation where all points of the subspace $A$ are equivalent to each other:

$$X / A = X / \sim, \ \ \ x \sim y, \ \text{if} \ \{x, y\} \subset A.$$

\textbf{Definition.}   \index{topological product}
The \textit{topological product} of spaces $X, Y$ is the Cartesian product $X\times Y$, where the topological structure consists of unions of products of open sets in $X$ and $Y$.

Let $X, Y$ be topological spaces with $X\cap Y=\emptyset$. On $X \cup Y$, a topology called \textit{disjoint union} is defined, where the open sets are unions of an open set in $X$ with an open set in $Y$.

\textbf{Definition.}  
\index{attachment}
If $A\subset Y$, and $f: A \to X$ is a continuous mapping, we say that $$X \cup_fY=(X\cup Y)/\sim, \ \ \text{where} \ \ x\sim f(x), \ x\in A$$
is obtained by \textit{gluing} (or \textit{attaching}) $Y$ to $X$ via the mapping $f$.

\textbf{Definition.}  
\index{wedge}
Let $x_0\in X, y_0\in Y$ be distinguished points in the spaces $X$ and $Y$. The \textit{wedge} of spaces is the disjoint union with the distinguished points glued together:

$$X \vee Y=(X\cup Y)/\sim, \ \ \text{where} \ \ x_0\sim y_0.$$

\textbf{Definition.}   \index{cylinder}
The cylinder over the space $X$ is defined as its topological product with the interval

$$\text{Cyl}(X)=X\times [0,1].$$

\textbf{Definition.}   \index{cone}
The cone over the space $X$ is the quotient space of the cylinder by the upper base

$$\text{Con}(X)=(X\times [0,1])/(X\times \{1\}).$$

\textbf{Definition.}   \index{inclusion}
An embedding of the space $X$ into the space $Y$ is a homeomorphism from $X$ to some subspace of $Y$.

\subsection{Classic Spaces and Their Notations}

We refer to the following as classic (standard) spaces:

1. \textit{The interval} $$I=[0,1];$$

2. \textit{The n-dimensional disk} $$D^n=\{(x_1,x_2, \ldots, x_n ) \in \mathbb{R}^n | \sum_{i=1}^n x_i^2 \le 1\};$$

3. \textit{The n-dimensional sphere} $$S^n=\{(x_1,x_2, \ldots, x_{n+1} ) \in \mathbb{R}^{n+1} | \sum_{i=1}^{n+1} x_i^2 = 1\};$$

In particular, the one-dimensional sphere is a circle $S^1$.

4. \textit{The n-dimensional real projective space}

$$\mathbb{R}P^n= (\mathbb{R}^{n+1} \setminus 0) / \sim,$$ $$ (x_1,x_2, \ldots, x_{n+1}) \sim (k x_1,k x_2, \ldots,k x_{n+1} ), k \in \mathbb{R};$$

5. \textit{The Möbius strip}

$$Mo=I^2/\sim, \ \ (a,0)\sim (1-a,1), \ a \in I=[0,1];$$

\begin{figure}[ht!]
\center{\includegraphics[width=0.55\linewidth ]{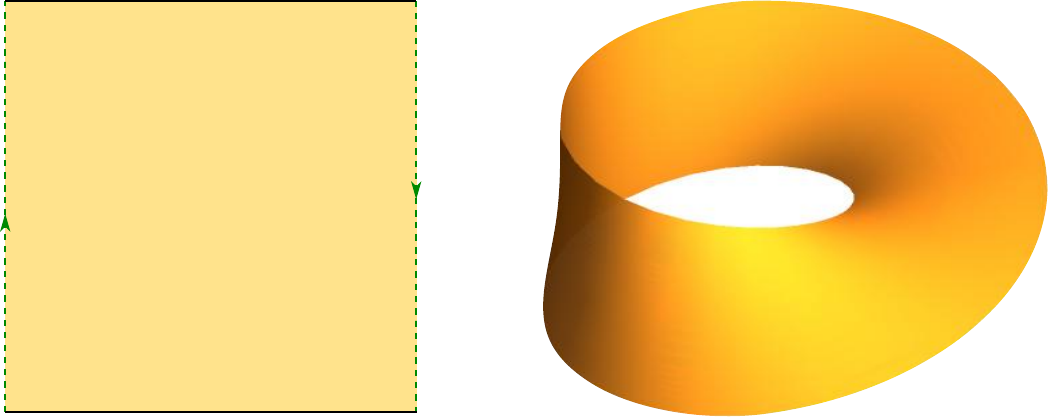}}
\caption{ The Möbius strip $Mo$}
\label{Mo}
\end{figure}

6. \textit{The torus} $$T^2=I^2/\sim, \ \ (a,0)\sim (a,1), \ (0,b) \sim (1,b) \ a,b \in I=[0,1];$$

\begin{figure}[ht!]
\center{\includegraphics[width=0.6\linewidth ]{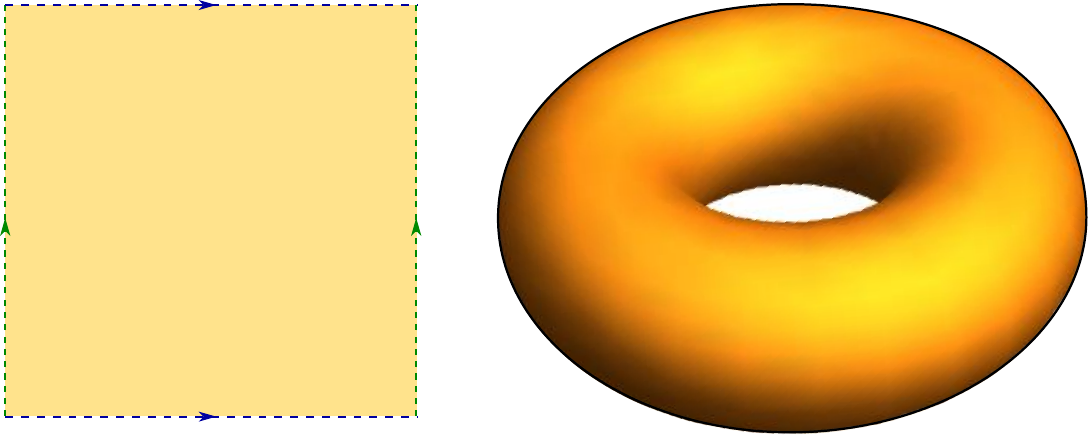}}
\caption{The two-dimensional torus $T^2$}
\label{t2}
\end{figure}

7. \textit{The Klein bottle}

$$Kl=I^2/\sim, \ \ (0,a) \sim (1, a), (b,0)\sim (1-b,1), \ a,b\in I.$$

\begin{figure}[ht!]
\center{\includegraphics[width=0.5\linewidth ]{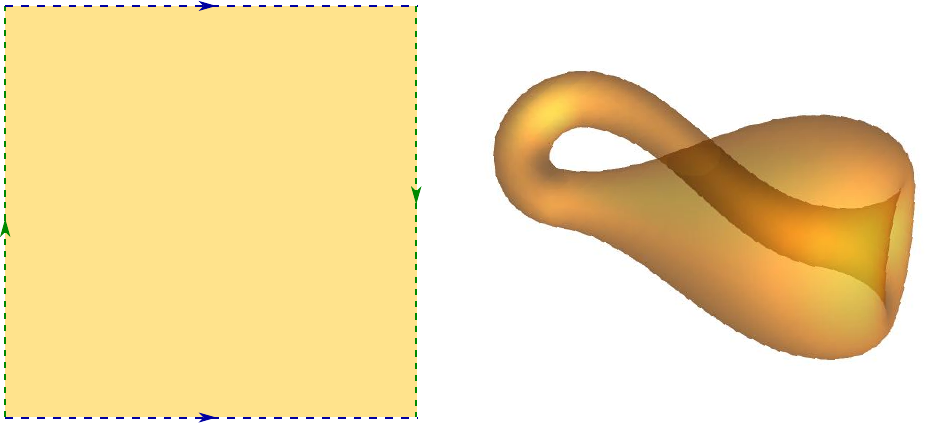}}
\caption{ The Klein bottle $Kl$}
\label{kl}
\end{figure}

\textbf{Problems. }

1. Prove that the Klein bottle \( Kl = I^2/\sim = [0,1]\times [0,1]/\sim \), where \( (0,a) \sim (1,a) \) and \( (a,0) \sim (1-a,1) \) for all \( a \in [0,1] \), is homeomorphic to the quotient space \( I^2/\sim = [0,1]\times [0,1]/\sim \), where \( (0,a) \sim (a,1) \) and \( (a,0) \sim (1,a) \).

2. Prove that if you seal the hole in a Klein bottle with a Möbius strip, you get a space that is homeomorphic to a torus with a hole sealed with a Möbius strip.

3. Show that if you glue opposite sides of a hexagon using homeomorphisms that are axial symmetries, you end up with a torus.

4. What happens if you cut the Möbius strip along the middle line?

\section{Smooth Manifolds and Their Mappings}

\textbf{Definition.}   \index{manifold}
A subset \( M \) of Euclidean space \( \mathbb{R}^N \) is called \textit{a submanifold} of dimension \( n \) if for every point \( p \in M \), there exists a neighborhood \( U \) in \( M \) and a homeomorphism \( h: U \to V \subset \mathbb{R}^n_+ \), where \( V \) is an open subset of the upper half-space $$ \mathbb{R}^n_+ = \{(x_1, x_2, \ldots , x_n) \in \mathbb{R}^n, x_n \ge 0\}.$$ The pair \( (U,h) \) is called a chart at point \( p \). By \textit{a manifold} we mean a topological space that is homeomorphic to a submanifold of Euclidean space.
Often, in notation, the upper index denotes the dimension of the manifold. The union of points with $$ h^{-1}(\{(x_1, x_2, \ldots , x_n) \in \mathbb{R}^n, x_n = 0\}) $$ forms \textit{the boundary} \( \partial M \) of the manifold \( M \). Points that do not belong to the boundary are called \textit{interior points} of the manifold.

\textbf{Definition.}   \index{атлас}
In the context of manifold theory, a set of charts that covers a manifold is called an \textit{atlas}. An atlas is called \textit{smooth} (of class $C^k$) if for any two charts $(U,h), (V,g)$ from it, the transition functions $g \circ h^{-1}: h(U\cap V) \to \mathbb{R}^n$ are smooth (of class $C^k$). A manifold that has a smooth atlas is referred to as \textit{a smooth manifold}.

\begin{figure}[ht!]
\center{\includegraphics[width=0.6\linewidth ]{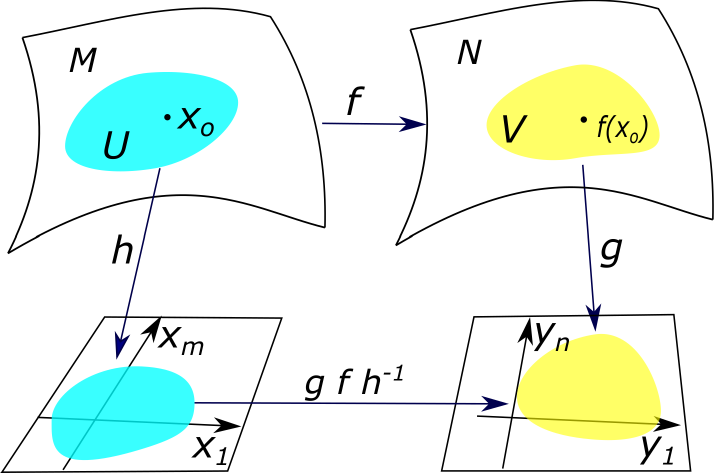}} 
\caption{ coordinate mapping record
 }
\label{fkoord}
\end{figure}

\textbf{Definition.}   \index{diffeomorphism}
A mapping $f:M^n \to N^k$ between smooth manifolds is termed a \textit{diffeomorphism} (of class $C^k$) if it is bijective, smooth, and the inverse mapping is also smooth (of class $C^k$).

A smooth function on a manifold $M^n$ is a smooth mapping $f:M^n \to \mathbb{R}$. 

\textbf{Definition.}   \index{submanifold}
A \textit{submanifold} of dimension $k$ in Euclidean space $\mathbb{R}^n$ is defined as a subset $M \subset \mathbb{R}^n$ such that for every point $x \in M$, there exists a neighborhood $U$ in $\mathbb{R}^n$ and a diffeomorphism $g:U \to V \subset \mathbb{R}^n$ such that $g(U \cap M) = g(U) \cap \mathbb{R}^k$, where $\mathbb{R}^k = \{(x_1, \ldots, x_n) \in \mathbb{R}^n: x_{k+1} = 0, x_n = 0\}$. The number $n-k$ is called the codimension of $M$. 
  
\begin{figure}[ht]
\center{\includegraphics[height=4.2cm]{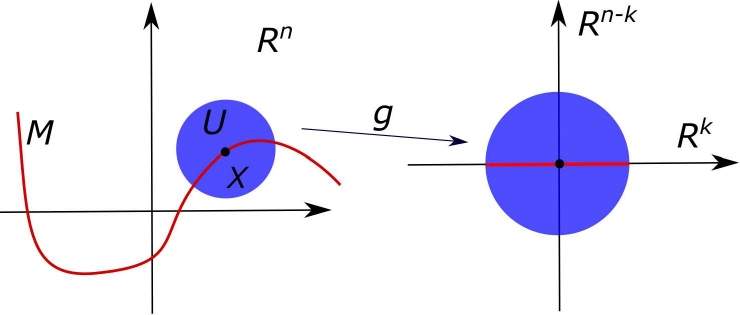}}
\caption{A mapping of a point neighbourhood}
\label{pic1}
\end{figure}

\begin{theorem}

(about the implicit function for submanifolds) Let a system of equations be given in the space $\mathbb{R}^n$

$$ f_i(x_1, \ldots, x_n) = 0, \ i=1,\ldots,m,
$$ 
where $f_i$ are smooth functions and $M$ is the set of solutions to this system. If the rank of the matrix $J = \left( \frac{\partial f_i}{\partial x_j} \right)$ is equal to $k$ everywhere on the set $M$, then $M$ is a submanifold of dimension $n-k$.
\end{theorem}
Unless otherwise stated, for submanifolds of Euclidean space, the atlas consists of sufficiently small neighborhoods of points and their projections onto the corresponding coordinate subspaces.

\textbf{Definition.}   \index{closed manifold}
A manifold is called \textit{closed} if it is compact and has an empty boundary.

\begin{theorem} (about the boundary) The boundary of a manifold of dimension $n$ is a manifold without boundary of dimension $n-1$.
\end{theorem}

\textbf{Corollary.} The boundary of a compact manifold is a closed manifold.

\begin{figure}[ht!]
\center{\includegraphics[width=0.6\linewidth ]{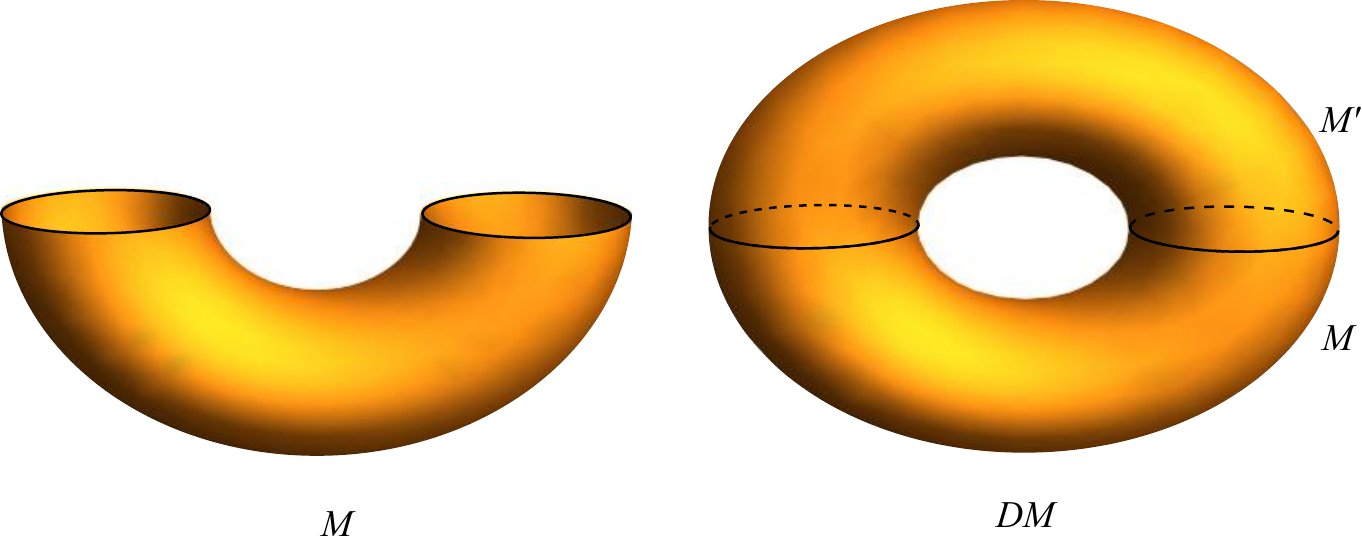}} 
\caption{ doubling of the manifold}
\label{DT2}
\end{figure}

\textbf{Definition.}   \index{doubling of a manifold}
The doubling $DM$ of a manifold $M$ is defined as the disconnected union of the manifold with itself, where the boundary points are glued together by an identity mapping:
$$DM = (M\cup M')/ \sim, \ \ \text{where} \ \ x\sim h(x), \ x\in\partial M, $$
here $h: M\to M'$ is some fixed homeomorphism (see Fig.\ref{DT2}).

\textbf{Definition.}   \index{curve} Let $M$ be a smooth manifold, $p \in M$. A \textit{curve} (starting at the point $p$) is defined as a smooth mapping $c:(a,b) \to M$, such that $a<0<b$, and $c(0)=p$.
  
Each such curve define the mapping $C^{\infty} (M)\to \mathbb{R}$ as a directional derivative at the point $p$.

\textbf{Problem.} 
Prove that this mapping is differentiation, i.e., it is a linear mapping and for any $f, g \in C(M): c(f g)= c(f) g + f c(g)$.

\textbf{Definition.}   Let $(U, h)$ be a chart at the point $p$ that provides local coordinates $x_1,\ldots , x_n$. Then the curve $c$ is represented by the usual curve $h(c(t))=\{x_1(t), \ldots , x_n(t)\}$ in $\mathbb{R}^n$. Two curves $c$ and $b$ are called equivalent if their velocity vectors (first derivatives) are equal in some chart.

\textbf{Problem.}  Prove that this definition does not depend on the choice of chart (Hint: the coordinates of the vectors are multiplied by a matrix of partial derivatives of the transition functions).

\textbf{Problem.}  Prove that curves are equivalent if and only if the induced differentiations coincide.

\textbf{Definition.}   \index{tangent vector} A \textit{tangent vector} is defined as the equivalence class of a curve under the equivalence relation introduced in the previous definition. The set of tangent vectors at the point $p$ is called the tangent space at the point $p$ and is denoted by $T_pM$.

Let $\partial_i$ be the tangent vector to the coordinate line $$x_1=0,\ldots, x_{i-1}=0, x_i=t, x_{i+1}=0,\ldots,x_n=0.$$ Coordinate-wise addition of curves and multiplication by a scalar induce linear operations on the space $T_pM$. Thus, the space $T_pM$ is an n-dimensional vector space with a basis $\{\partial_1,\ldots, \partial_n\}$. If a curve $c$ is given by the equation $x(t)= \{x_1(t), \ldots, x_n(t)\}$, then the coordinates of the tangent vector to it in this basis  is $$\{x_1'(t), \ldots, x_n'(t)\}.$$

\begin{figure}[ht!]
\center{\includegraphics[width=0.6\linewidth ]{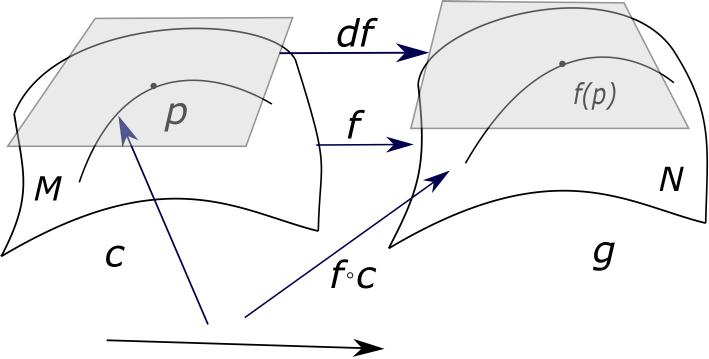}} 
\caption{differential}
\label{dgt8}
\end{figure}

\textbf{Definition.}   \index{differential of a mapping}
Let $f: M\to N$ be a smooth mapping of smooth manifolds, with $p\in M$. Then the mapping $f$ induces a mapping $$df_p: T_pM \to T_{f(p)}N$$ of tangent spaces: if $c$ is a curve that defines a vector at the point $p$, then the composition $f\circ c$ is a curve that defines a vector at the point $f(p)$. The mapping $df_p$ is called the \textit{differential} of the mapping $f$ at the point $p$. Sometimes the mapping $df_p$ is referred to as the tangent mapping and denoted $T_pf$ (see Fig. \ref{dgt8}).

\textbf{Problem.} 
Prove that the mapping $d f_p$ is defined correctly, is a linear mapping, and in local coordinates can be expressed by the formula $df_p(\{a_1, \ldots , a_n\})=$, where $y_i=f_i(x_1,\ldots , x_n)$ is the local representation of the mapping $f$.

\textbf{Problem.} Prove that if $M_1, M_2, M_3$ are smooth manifolds, and $f: M_1\to M_2$, $g: M_2\to M_3$ are smooth mappings, then $d (g \circ f)_p= d g_{f(p)} \circ d f_p$.

From the last two problems, the following follows.

\textbf{Proposition.} If $f$ is a diffeomorphism, then $d f_p$ is an isomorphism.

\textit{Definition.} A point $p\in M$ is called a regular point of the mapping $f: M^n \to N^k$ if the rank of the matrix $\left( \frac{\partial f_i}{\partial x_j} \right)$ is equal to $k$. Irregular points are called \textit{critical} or \textit{singular}.

\textbf{Definition.} A mapping $f: M \to N$ is called an \textit{ immersion} if for any point $p \in M$, the differential of the mapping $d f: T_pM \to T_{f(p)}N$ is a monomorphism, i.e., an isomorphism onto a subspace of the space $T_{f(p)}N$. A mapping $f: M \to N$ is called an \textit{embedding} if it is an immersion and an injective mapping.

\textbf{Definition.} \textit{The tangent bundle} of a manifold $M$ is a non-connected union of tangent spaces along with the canonical projection $\pi : TM\to M$ such that $\pi (T_p M)=\{p\}$. The tangent bundle $TM$ is equipped with the structure of a manifold: each chart $(U, h)$ on $M$ induces a chart $(U \times \mathbb{R}^n, h)$, where the homeomorphism $H$ is given by the formula $$H (p, v) := (h(p), (v_1,\ldots , v_n)).$$ Here, $v_i$ are the coefficients (coordinates) of the vector $v$ in its decomposition with respect to the basis $\{ \partial x_1,\ldots , \partial x_1\}$ of the space $T_pM$.

\textbf{Problem.} Show that $TM$ is a smooth manifold (of dimension $2n$), and that the mapping $\pi$ is a smooth mapping of manifolds.

Let the manifold \( M \) be embedded in \( \mathbb{R}^N \) (the existence of such an embedding follows from the Whitney embedding theorem). Then the tangent bundle \( TM = \{(p,v) | p \in M, v \in T_pM\} \) is a manifold embedded in \( \mathbb{R}^{2N} \).

\textbf{Definition.}   \index{transversal} Let \( f: M \to N \) be a smooth map between smooth manifolds, and let \( L \) be \textit{a submanifold }of the manifold \( N \). The map \( f \) is called \textit{transversal} to \( L \) at the point \( p \in M \) if \( f(p) \notin L \) or \( df(T_pM) + T_{f(p)}L = T_{f(p)}N \). The map \( f \) is \textit{transversal} to \( L \) if it is transversal at every point \( p \in M \). A submanifold \( L_1 \subset N \) is called \textit{transversal} to the submanifold \( L_2 \subset N \) if the embedding map \( i: L_1 \to N \) is transversal to \( L_2 \).
  
It is obvious that if \( L_1 \) is transversal to \( L_2 \), then \( L_2 \) is transversal to \( L_1 \). Transversal submanifolds are also said to be in general position.

\textbf{Example.} 
If \( \text{dim} L_1 + \text{dim} L_2 < \text{dim} N \), then \( L_1 \) is transversal to \( L_2 \) if and only if \( L_1 \cap L_2 = \emptyset \). If \( f: M \to N \) is a submersion, then \( f \) is transversal to any submanifold \( L \subset N \). In a two-dimensional manifold (surface), two curves are transversal if they have no points of tangency.

\textbf{Problem.} 
Which of the following subsets of Euclidean space are (smooth) manifolds (prove your answer):

1) a circle; 2) the union of two circles that have one common point; 3) two circles that have two common points; 4) two non-intersecting circles; 5) a triangle; 6) a tetrahedron; 7) two triangles that share a common edge; 8) the union of a sphere and a circle that do not intersect?

\section{Vector Fields and Dynamical Systems}

\textbf{Definition.}   \index{vector field} We say that a \textit{vector field} $X$ is defined on a manifold $M$ if at each point $p$ there is a tangent vector $v_p \in T_pM$. Thus, a \textit{vector field} is a mapping $X: M \to TM$ such that the composition $\pi \circ X = Id$ is the identity mapping.

In each coordinate chart with coordinates $x_1, \ldots, x_n$, the vector field is defined by the functions $X_i = X_i(x_1, \ldots, x_n)$ — the components of the vectors with respect to the basis $\{\partial x_1, \ldots, \partial x_n\}$. A vector field is said to be smooth if all functions $X_i$ are smooth. Henceforth, unless stated otherwise, we will assume all vector fields are smooth. Since each vector can be viewed as differentiation (derivative in the direction), the task of defining a vector field is equivalent to the task of differentiating smooth functions $X: C^{\infty}(M) \to C^{\infty}(M)$.

\textbf{Definition.}   \index{trajectory} A trajectory of a vector field is a curve such that the tangent vector at any point on this curve coincides with the vector field $X$ at that point.

In a coordinate chart with coordinates $x_1, \ldots, x_n$, a trajectory is a solution to the system of ordinary differential equations:
$$x_1' = X_1(x_1, \ldots, x_n) $$
$$ \dots \dots$$
$$x_n' = X_n(x_1, \ldots, x_n).$$

\begin{figure}[ht!]
\center{\includegraphics[width=0.7\linewidth ]{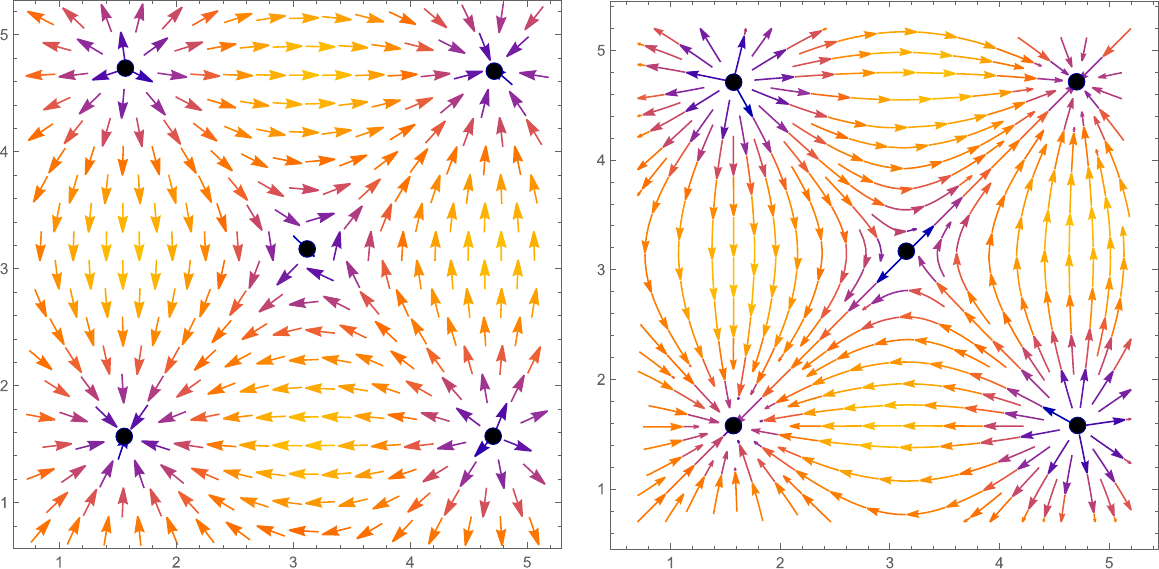}}
\caption{vector field $\{\sin u \cos v, \cos u \sin v \}$ and its trajectories}
\label{vecftr}
\end{figure}
From the theorem on the existence and uniqueness of solutions to the system of differential equations, it follows that through each point there passes a unique trajectory $x_i = x_i(t)$.

\textbf{Definition.}   \index{complete trajectory} A trajectory is called complete if it is defined for all $t \in \mathbb{R}$. A vector field is called complete if through each point of the manifold there passes a complete trajectory.
  
Every vector field on a closed manifold is complete.

\textbf{Definition.}   Possible three types of trajectories on a manifold:

1) \index{singular point} \textit{singular points} -- a continuous mapping of a line into a point. These are the points where the vector field equals 0;

2) \index{simple trajectory} \textit{simple trajectories} -- embeddings of a line or an interval;

3) \index{closed trajectory} \textit{closed (or periodic) trajectories} -- closed curves.

It should be noted that simple trajectories can either be homeomorphic to a line (for example, in the plane for the field $X=\{1,0\}$ the trajectories are horizontal lines) or not (for example, on a torus with parametrization $$\{(2+ \cos u)\cos v, (2+ \cos u)\sin v, \sin u\}, u \in [0, 2 \pi], v \in [0, 2 \pi],$$
each trajectory of the field $\{1, a\}$, which is defined in the local coordinate system $u,v$, for irrational $a$ is everywhere dense on the torus).

Let $M$ be a closed manifold. We denote by $\gamma_t(p), t \in \mathbb{R}$ -- such a trajectory that $\gamma_0(p)=p$. Fix $t\in \mathbb{R}$. Then we have a mapping $\gamma_t:M \to M$. We will show that this mapping is a diffeomorphism. Indeed, there exists an inverse mapping $(\gamma_t)^{-1}=\gamma_{-t}$. The smoothness of the mapping $\gamma_t$, as well as its inverse, follows from the fact that the solution of the differential equation depends smoothly on the initial conditions.

\textbf{Definition.}   \index{} A set of diffeomorphisms $\gamma_t:M \to M, t \in \mathbb{R}$ is called a one-parameter group of diffeomorphisms or a flow if

1) $\gamma_0=Id$,

2) $\gamma_{-t}=(\gamma_t)^{-1}, t \in \mathbb{R}$,

3) $\gamma_{t+s}(p)=(\gamma_t(\gamma_s(p)), t, s \in \mathbb{R}, p \in M$.

It is obvious that every vector field on a closed manifold generates a flow. Conversely, every flow defines curves $\gamma_t(p)$ for each $p$, and therefore, a vector field consisting of tangent vectors to these curves. Thus, a vector field can be defined in the following ways:

1) A mapping $X: M \to TM$ such that $\pi \circ X = Id$,

2) The differentiation of smooth functions $X: C^{\infty}(M) \to C^{\infty}(M)$,

3) A system of differential equations (*),

4) A flow (a one-parameter group of diffeomorphisms) $\gamma_t: M \to M, t \in \mathbb{R}$.

The set of all smooth vector fields defined on the manifold $M$ will be denoted as $\Gamma(M)$.

\textbf{Problem.} 
Let $p$ be a non-singular point of the vector field $X$ $(X_p \ne 0)$. Prove that there exists a coordinate system $(x_1, \ldots, x_n)$ in the point $p$ such that $X = \{1, 0, \ldots, 0\}$.

\textbf{Definition.}   \index{Riemannian metric} If at each point $p \in M$ a scalar product (a bilinear symmetric positive definite mapping) is defined in the tangent space $T_pM$ that continuously depends on the point, then it is said that a \textit{Riemannian metric} is defined on the manifold.
  
In a fixed chart $(U, h)$, the Riemannian metric is given by the matrix $$(g_{ij}(p))_{i,j=1}^n, p \in M.$$
An example of a Riemannian metric is the first quadratic form of a surface in three-dimensional space.

\textbf{Definition.}   From the course of linear algebra, it is known that the scalar product $g:V \times V \to\mathbb{R}$ on the vector space $V$ defines its isomorphism to the dual space $V^*$: $$ v \in V \to v^* \in V^*, v^*(w)=g(v,w).$$ A vector $v$ is called \textit{dual} to the linear mapping $v^*$.
  
\textbf{Definition.}   \index{gradient field} If $f:M \to \mathbb{R}$ is a smooth function, then we define the vector grad $f(x) \in T_pM$ for $p \in M$ as dual to the linear mapping $df: T_pM \to \mathbb{R}$. The resulting vector field grad $f$ is called the \textit{gradient field} of the function $f$. This field depends on the choice of Riemannian metric.

If $M$ is an open subset of Euclidean space $\mathbb{R}^n$ with the standard metric, then $$ \text{grad} f =\nabla f =\left\{ \frac{\partial f}{\partial x_1}, \ldots , \frac{\partial f}{\partial x_n} \right\}.$$

\begin{figure}[ht!]
\center{\includegraphics[width=0.4\linewidth ]{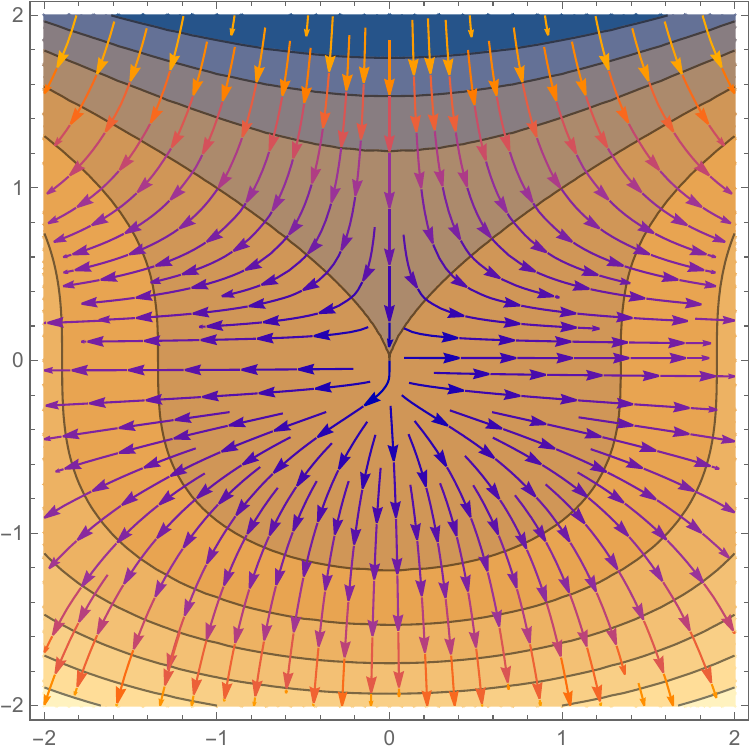}} 
\caption{level lines and the gradient vector field of the function $f(x,y)=x^2-y^3$.
 }
\label{grad_x2-y^3}
\end{figure}

For example, on the plane $\nabla (x^2-y^3)=\{2x,-3y^2\}$ (see \ref{grad_x2-y^3}).

Clearly, grad $f(p)=0$ if and only if the point $p$ is a critical point of the function $f$. At regular points, the vector grad $f(p)$ is orthogonal to the level surface (line) $f^{-1}(f(p))$.

\textbf{Problem.} 
 Prove that the function $f$ is non-decreasing along each trajectory of the vector field grad $f$.

\textbf{Definition.}   \index{non-degenerate point} A special point $p$ of a vector field $X$ 
is called \textit{non-degenerate} if the matrix $\left(\frac{\partial X_i}{\partial x_j}\right)_{i,j=1}^n$ does not have any eigenvalues whose real part is equal to 0.

\textbf{Definition.}   \index{stable manifold} A \textit{stable} \textit{manifold} at the point $p$ is defined as the set
$$S(p)=\{x \in M: \lim_{t\to \infty} \gamma_t(x)=p\}.$$ A \textit{unstable} manifold is the set $$U(p)=\{x \in M: \lim_{t\to - \infty} \gamma_t(x)=p\}.
$$   


\textbf{Definition.}   \index{Morse field} A \textit{gradient-like Morse-Smale field} (or \textit{Morse field}) is called a vector field $X$ such that

1) the set of critical points $\Omega (X)$ is finite and all critical points are non-degenerate,

2) for every point $x \in M:$ $$\lim_{t \to +\infty } \gamma_t(x) \in \Omega (X) \ \text{and} \ \lim_{t \to -\infty} \gamma_t(x) \in \Omega (X)$$ (each trajectory starts and ends at a critical point),

3) stable and unstable manifolds of critical points intersect transversally.

\textbf{Problem.} 
Find the gradient fields of the projections onto the coordinate axes of the following surfaces: 1) sphere, 2) torus, 3) hyperboloids of revolution, 4) paraboloids.

\section{Morse Functions}

Let $M$ be a smooth manifold, and let $f: M\to \mathbb{R}$ be a smooth function.

\textbf{Definition.}   \index{critical point} A point $p\in M$ is called a critical point of the function $f$ if the differential $df(p)=0$. The value $f(p)$ is called the critical value of the function $f$.  

If $(U, h)$ is a chart at the point $p$, $h(p)=0$, and $(x_1,\ldots , x_n)$ are local coordinates, then the point $p$ will be critical if
$$ \frac{\partial f}{\partial x_1}=\frac{\partial f}{\partial x_2}=\ldots =\frac{\partial f}{\partial x_n}=0.$$
\textbf{Definition.}   \index{non-degenerate point} A critical point $p$ of the function $f$ is called \textit{non-degenerate} if at it the Hessian matrix $H$ of second derivatives is non-degenerate:
$$\text{det} \left( \frac{\partial ^2f}{\partial x_i\partial x_j} \right) _{i,j=1}^n \ne 0.$$

%




\textbf{Problem.} 
Check that this definition does not depend on the choice of the local coordinate system at the point $p$.
 
\textbf{Definition.}   \index{Morse function} A smooth function is called a \textit{Morse function} if all its critical points are non-degenerate.

\textbf{Example.} 
The function $f: \mathbb{R}^2\to\mathbb{R}, f(x, y)=xy$ has one critical point $(0,0)$, which is non-degenerate. The function $f: \mathbb{R} \to \mathbb{R}, f(x)=x^3$ has one degenerate critical point at 0.
 
\begin{theorem} \textbf{(Morse).} Let $p$ be a non-degenerate critical point of a smooth function $f: M\to \mathbb{R}$. Then there exists a local coordinate system $x_1,\ldots , x_n$ at the point $p$, in which $p=(0,0,\ldots ,0)$ and
$$ f(x_1,\ldots , x_n)=f(p)-x_1^2-x_2^2-\ldots -x_{\lambda}^2+x_{\lambda+1} ^2+\ldots +x_n^2.$$
\end{theorem}
\textbf{Proof.}
\textit{Necessity.} Without loss of generality, let us assume that $p=0$ and $f(p)=f(0)=0$. Then for some neighborhood of $0$

$$f(x_1,\ldots, x_n)= \int_0^1\frac{df(tx_1, \ldots,tx_n)}{dt} dt=$$ $$=\int_0^1\sum_{i=1}^n \frac{df}{dx_i}(tx_1, \ldots,tx_n) x_i dt
= \sum_{i=1}^n x_i g_i(x_1,\ldots, x_n),$$  $\text{where} \ g_i(x_1,\ldots, x_n)=\int_0^1 \frac{df}{dx_i}(tx_1, \ldots,tx_n) dt.$

Since $0$ is a critical point of $f$, we have $g_i(0)=0$. Then
$$g_i(x_1,\ldots, x_n)= \int_0^1\sum_{j=1}^n \frac{dg_i}{dx_j}(sx_1, \ldots,sx_n) x_j ds=$$
$$= \sum_{j=1}^n x_j h_{ij}(x_1,\ldots, x_n),$$ where $h_{ij}(x_1,\ldots, x_n)=\int_0^1 \frac{dg_i}{dx_j}(sx_1, \ldots,sx_n) ds.$
Thus,
$$f(x_1,\ldots, x_n)= \sum_{i,j=1}^n x_i x_j h_{ij} (x_1,\ldots, x_n).$$

The matrix $(h_{ij}(0))$ is equal to the matrix $H(0)$, hence it is non-degenerate.
Using linear coordinate transformations, we will reduce the matrix $(h_{ij}(0))$
to diagonal form. The same transformations will be applied to the expression $$\sum_{i,j=1}^n x_i x_j h_{ij} (x_1,\ldots, x_n).$$ By completing the square, we will reduce it to a diagonal matrix and denote the expressions under the squares with new variables. Then, using the implicit function theorem, we prove the necessary conditions of the theorem.

\textit{Sufficiency.} All first partial derivatives of the function
$$ f(x_1,\ldots , x_n)=f(p)-x_1^2-x_2^2-\ldots -x_{\lambda}^2+x_{\lambda+1} ^2+\ldots +x_n^2 $$
are equal to 0 at the point $p$,
and the matrix $H$ is diagonal with numbers $\pm 2$ on the diagonal.
Therefore, $p$ is a non-degenerate critical point.

\textbf{Definition.}   \index{Morse index} The \textit{Morse index} of a non-degenerate critical point is the number $\lambda$ from the previous theorem. It does not depend on the choice of the coordinate system.

If the Morse index is equal to 0, then the point is a local minimum, and if it is equal to $n$, then it is a local maximum.

\textbf{Corollary.} All non-degenerate critical points are isolated.
\begin{figure}[ht!]
\center{\includegraphics[width=0.95\linewidth ]{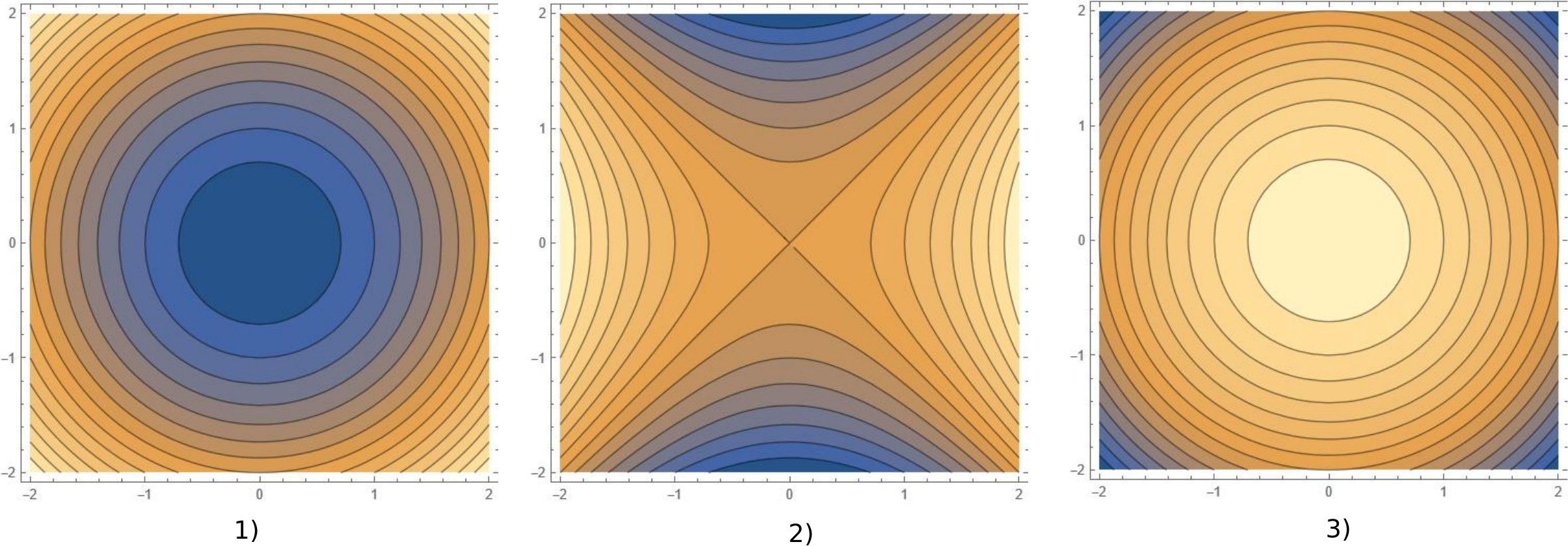}} 
\caption{nondegenerate critical points on 2-manifolds
 }
\label{Morse2}
\end{figure}

On two-dimensional manifolds, depending on the Morse index, critical points can be of three types (see Fig. \ref{Morse2}):

1) $\text{ind} = 0$, $f(x,y)=x^2+y^2$ -- local minimum, $\text{grad} f=\{2x,2y\}$ -- source;

2) $\text{ind} = 1$, $f(x,y)=x^2-y^2$ -- saddle, $\text{grad} f=\{2x,-2y\}$ -- saddle;

3) $\text{ind} = 0$, $f(x,y)=-x^2-y^2$ -- local maximum, $\text{grad} f=\{2x,2y\}$ -- sink.

\begin{figure}[ht!]
\center{\includegraphics[width=0.95\linewidth ]{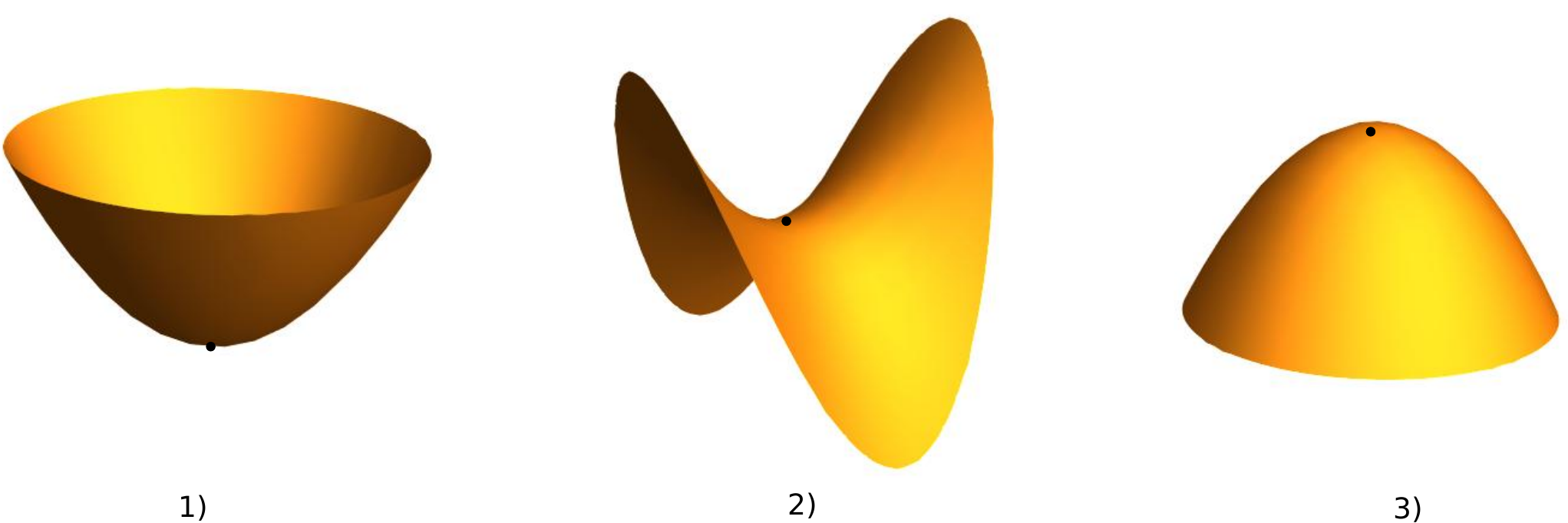}}
\caption{Graphs of the function in the vicinity of non-degenerate critical points on 2-manifolds}
\label{Mfunc-h}
\end{figure}

The graphs $z=f(x,y)$ for these three functions are depicted in Fig. \ref{Mfunc-h}.

\begin{figure}[ht!]
\center{\includegraphics[width=0.95\linewidth ]{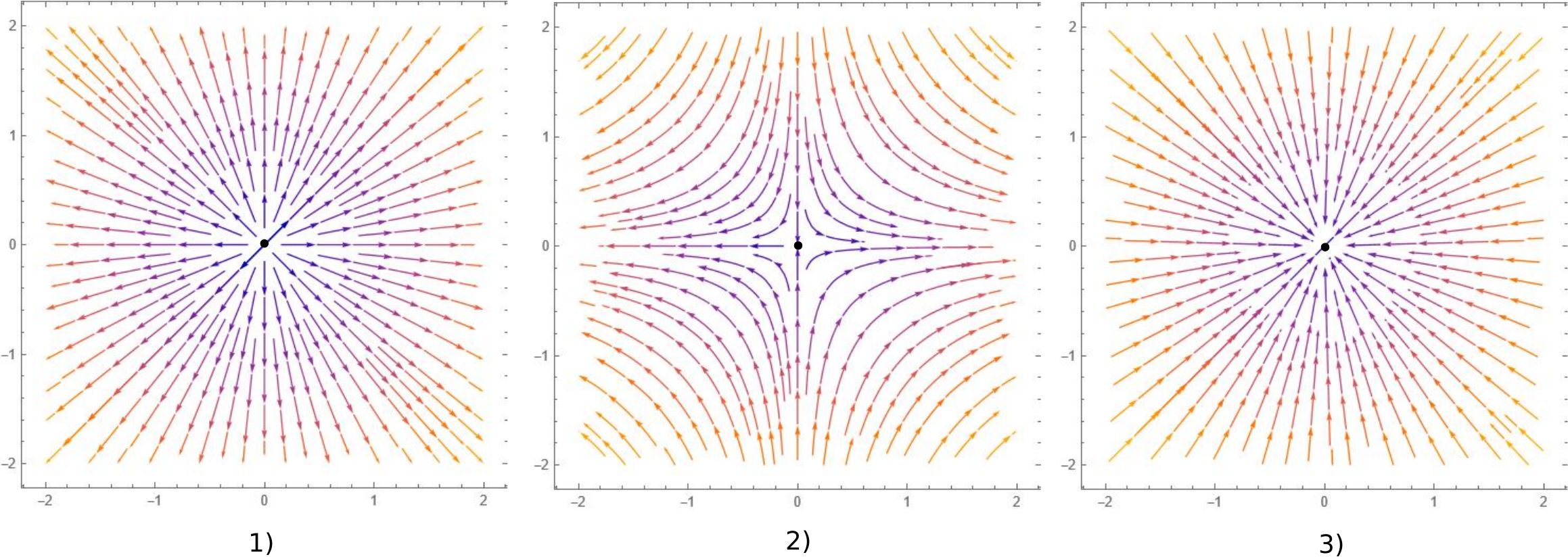}}
\caption{Trajectories of the gradient field of the Morse function in the vicinity of critical points}
\label{Mgrad}
\end{figure}

The gradient field of the Morse function on 2-manifolds in the neighborhoods of critical points is illustrated in Fig. \ref{Mgrad}.

\begin{figure}[ht!]
\center{\includegraphics[width=0.95\linewidth ]{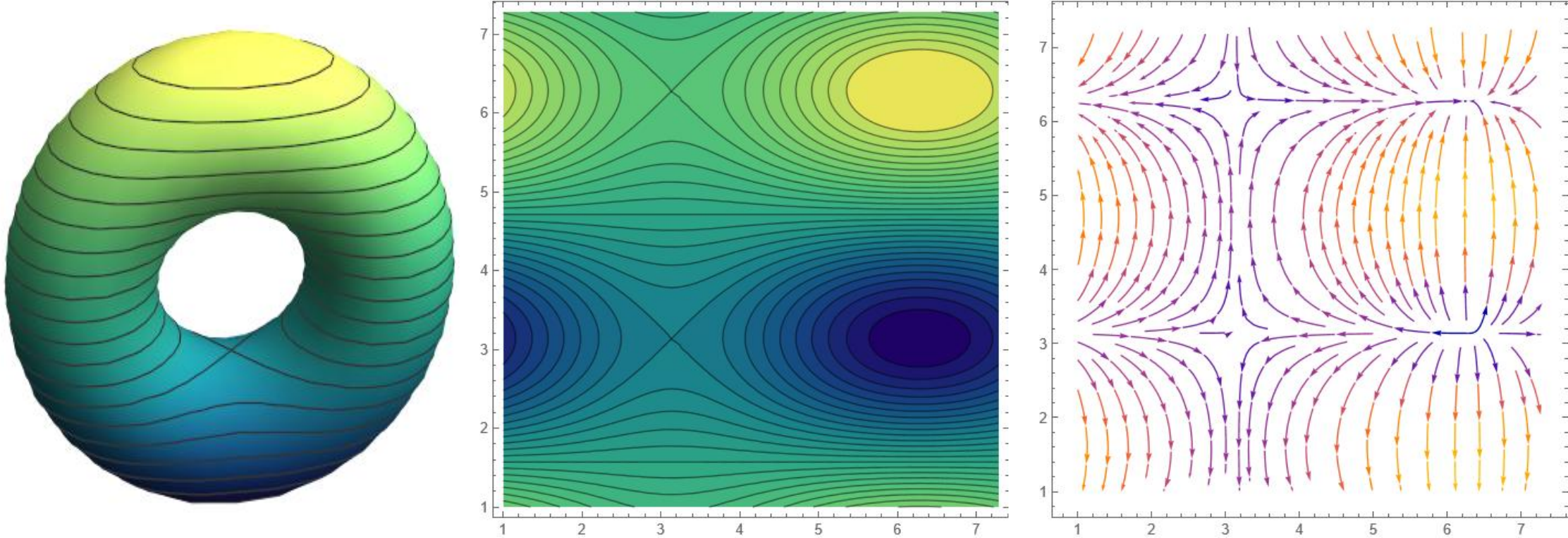}}
\caption{Level curves of the height function on the torus in three-dimensional space, along with the unfolding and trajectories of its gradient field}
\label{h-fun-t2}
\end{figure}

\textbf{Example.} Let the torus be defined by the parametric equation
$$\overline{r(u,v)}= \{(\cos (u)+2) \sin (v),\sin (u),(\cos (u)+2) \cos (v)\},$$
$u \in [1,1+2\pi], v\in [1,1+2\pi]$. Consider the height function on the torus $h(x,y,z)=z$. In local coordinates, it can be expressed as $h(u,v)= (\cos (u)+2) \cos (v)$. Its level curves and the trajectories of its gradient field are depicted in Fig. \ref{h-fun-t2}.

Let's find the critical points of the function $h$:
$$h'_v= -\sin v (2+\cos u)=0, \ \ v_1=0, v_2=\pi;$$
$$h'_u= -\sin u \cos v =0, \ \ u_1=0, u_2=\pi.$$
Thus, we have 4 critical points: $$p_1=(0,0), p_2=(0, \pi), p_3=(\pi, 0), p_4=(\pi,\pi).$$

Let's check the non-degeneracy of the critical points:

$$ H_{11}=-\cos u \cos v, H_{12}=0, H_{22}=-\cos v (2+ \cos u).$$ 

$$H(p_1)=
\begin{pmatrix}
  -1& 0\\
  0& -3
\end{pmatrix}, \ \ \ \ \ \
H(p_2)=
\begin{pmatrix}
  1& 0\\
  0& 3
\end{pmatrix},$$
$$ H(p_3)=
\begin{pmatrix}
  1& 0\\
  0& -1
\end{pmatrix}, \ \ \ \ \ \
H(p_4)=
\begin{pmatrix}
  -1& 0\\
  0& 1
\end{pmatrix},
$$

Conclusion: all critical points are non-degenerate, $h$ is a Morse function, $p_1$ is a maximum point (ind=2), $p_2$ is a minimum point (ind=0), and $p_3$ and $p_4$ are saddle points (ind=1).

Let us find the gradient field of this function in the Riemannian metric $ds^2=du^2+dv^2$ (see Fig. \ref{h-fun-t2}):
$$\text{grad } h=\nabla h= \{ -\sin u \cos v, -\sin v (2+\cos u)\}.$$
From the non-degeneracy of the critical points of the function, it follows the non-degeneracy of the singular points of the gradient field. Since the function increases along each trajectory, the boundary sets of the trajectories are the singular points. We will check the transversality condition of stable and unstable manifolds.

For the singular point $p_3=(\pi, 0)$, in its neighborhood, the gradient field is close (after Taylor expansion) to the field $\{ u, -v\}$. Therefore, the stable manifold is given by the equation $u=0$ (at each such point, the first coordinate of the gradient vector equals 0).

For the singular point $p_4=(\pi,\pi)$, in its neighborhood, the gradient field is close (after Taylor expansion) to the field $\{- u, v\}$. Therefore, the unstable manifold is also given by the equation $u=0$ (at each such point, the first coordinate of the gradient vector equals 0). Thus, the stable manifold $p_3$ intersects with the unstable manifold $p_4$ along two trajectories (saddle connections). This means that the gradient field in this Riemannian metric is not a Morse-Smale field.

Let us denote $M_y=f^{-1}((-\infty, y])$, $\partial M_y=N_y=f^{-1}(y)$.
\begin{lema} If the interval $[y, z]$ does not contain critical values, then $M_y$ is homeomorphic to $M_z$, and $N_y$ is homeomorphic to $N_z$.
\end{lema}
\textbf{Proof.}
Consider the vector field grad $f$. Since $L=f^{-1}([y, z])$ does not contain critical points and the function increases along the trajectories, every trajectory passing through $L$ intersects each of $N_y$, $N_z$ at one point. Hence, we have a bijection between $N_y$, $N_z$. From the theory of differential equations, it follows that this bijection is continuous in both directions, meaning it is a homeomorphism.

A homeomorphism between $M_y$ and $M_z$ can be constructed by stretching each segment of the trajectory

$\gamma(x) \cap f^{-1}([y-\varepsilon, y])$ to the segment of the trajectory $\gamma(x)\cap f^{-1}([y-\varepsilon, z])$.


As a consequence, we have:

\begin{theorem}(Reeb) If a Morse function on a closed manifold $M^n$ has two critical points, then this manifold is homeomorphic to the sphere $S^n$.

\end{theorem}

\textbf{Proof.}
Let the function $f:M^n\to \mathbb{R}$ have two critical points, $p$ -- the point of minimum, and $q$ -- the point of maximum, $f(p)=a$, $f(q)=b$. Consider a coordinate system at the point $p$, in which the function has a canonical form. Then there exists $c$, $a<c,b$, such that $f^{-1}(c)$ is homeomorphic to $S^{n-1}$, ($g:f^{-1}(c) \to S^{n-1}$ -- this is the homeomorphism). On the sphere $S^n$, consider the height function $h:S^n \to \mathbb{R}$. Then $h^{-1}(0)=S^{n-1}$. Construct a map
$M^n\to S^n$, that maps the trajectories of the gradient field of the function $f$ to the trajectories of the gradient field of the function $h$. Here, the source maps to the source, and the sink to the sink. The bijection between other trajectories is determined by the mapping $g$. Choose as the initial point (the point with parameter $t=0$) on the trajectories on the manifold $M^n$ points from the level $f^{-1}(c)$, and on the sphere $S^n$ the initial points lie on $S^{n-1}$. Then the bijection between the trajectories together with the condition of parameter preservation defines a bijection between $M^n$ and $S^n$. By the theorem on the continuous dependence of the solution of a differential equation on initial conditions, the constructed bijection is a continuous map. Similarly, the inverse map is continuous. Thus, we have constructed the desired homeomorphism.


\textbf{Problem.} 
Prove that a Morse function on a compact manifold has a finite number of critical points.

\textbf{Problem.} 
Which of the following functions are Morse functions:

1) projections of the ellipsoid (hyperboloids, paraboloids, torus) onto the coordinate axes; 2) the distance function (square of the distance) to the origin on these surfaces; 3) the distance function to the coordinate axis on these surfaces,
4) The function \(f(x:y:z) = \frac{x^2+2y^2+3z^2}{x^2+y^2+z^2}\) in homogeneous coordinates on the projective plane, 5) The function \(f(x:y:z) = \frac{x^2-y^2}{x^2+y^2+z^2}\) in homogeneous coordinates on the projective plane, 6) Functions \(\sin(u)\cos(v)\sin(w)\) and \(\sin(u) \tan(v) \cos(w)\) in local coordinates on a three-dimensional torus? Find the critical points, and for non-degenerate critical points -- their Morse indices.

\section{Cell Complexes}

\textbf{Definition.}   \index{cell complex} A space \(K\) is called a cell complex (or CW complex, or cellular complex) if for any two points in \(K\), there exist neighborhoods that do not intersect and a partition of \(K\) into cells is given: so that for each cell \(e_i^k\) (where \(k\) is the dimension of the cell) there exists a continuous mapping of a closed \(k\)-dimensional disk \(D^k\) into \(K\), whose restriction to the interior of \(D^k\) is a homeomorphism onto \(e_i^k\), and such that

C) the boundary of each cell is contained in a finite union of cells of lower dimensions (closure finite),

W) a set \(A\) is closed in \(K\) if and only if its intersection with each cell is closed in the topology of these cells (weak topology).

\textbf{Definition.}   \index{dimension of the complex} The dimension of a cell complex is the highest dimension of its cells. The union of all cells in the space \(K\) whose dimension does not exceed \(n\) is called the n-skeleton and is denoted \(K^n\).

Every cell complex can be built in such an inductive manner:

1) \(K^0\) is a discrete space, where each point is a 0-cell;

2) \(K^n\) is obtained from \(K^{n-1}\) by attaching a disjoint union of \(n\)-dimensional disks \(D^n\) via continuous functions (attachment maps) \(f_i: \partial D^n \to K^{n-1}\);

3) \(K=\cup_n K^n\) (with the weakest topology as in W).

Examples of cell complexes.

\textbf{Example.} 
Sphere $S^n$. The most typical cell decompositions of the sphere are

a) $S^n=e^0\cup e^n, f(\partial D^n)=e^0$;

b) $S^n=e_1^0 \cup e_2^0\cup e_1^1\cup e_2^1\cup \ldots \cup e_1^n\cup e_2^n, e_1^k=\{x=(x_1,x_2,\ldots ,x_{n+1}) \in S^n | x_{k+1}>0, x_{k+2}=\ldots =x_{n+1}=0\}, e_2^k=\{x=(x_1,x_2,\ldots ,x_{n+1}) \in S^n | x_{k+1}<0, x_{k+2}=\ldots =x_{n+1}=0\}, f_i^k:\partial D^k \to S^{k-1} $ -- homeomorphism;

c) $S^n=\partial I^{n+1}$ -- the boundary of an $n+1$-dimensional cube, where the cells are vertices, edges, and faces.

\textbf{Example.} 
Disk $D^n=e^n \cup S^{n-1}$, where $S^{n-1}$ has one of the cell decompositions described above.

\textbf{Example.} 
Projective space $\mathbb{R} P^n=S^n/\sim$, where $x\sim -x$. From decomposition b) of the sphere, after identifying diametrically opposite points, we get $\mathbb{R} P^n=e^0\cup e^1\cup \ldots \cup e^n$, where $e^k= (e_1^k\cup e_2^k)/\sim, x\sim -x$.

\textbf{Example.}  Each compact surface can be obtained from a polygon by gluing pairs of corresponding sides. 0-cells will be the vertices of the polygon (after gluing), 1-cells will be the edges, and the 2-cell will be the interior of the polygon.

Let's describe the standard gluing of a polygon, where all vertices are glued into one. For a closed oriented surface in a quadrilateral with 4n sides, the sides 1 and 3, 2 and 4, 5 and 7, 6 and 8, etc., are glued together by homeomorphisms that reverse boundary orientation. For a closed non-oriented surface in a digon, adjacent sides (1 and 2, 3 and 4, etc.) are glued together by homeomorphisms that preserve boundary orientation.

\textbf{Definition.}   \index{cellular subspace}
A subset $A$ of a cell space $K$ is called a \textit{cell subspace} if $A$ is a cell space and every cell in $A$ is a cell of space $K$.

\section{Homotopic Mappings}

Let's equip the set $C(X,Y)$ of all continuous mappings from a topological space $X$ to a topological space $Y$ with the compact-open topology. The open sets of this topology can be obtained through finite intersections and arbitrary intersections of sets $$U(A,B)=\{f\in C(X,Y)|f(A)\subset B, A - \text{compact}, B - \text{open}\}.$$ For locally compact Hausdorff spaces $X, Y, Z$ (these are primarily the spaces considered), the natural mapping $C(Z, C(X,Y)) \to C(X\times Z, Y)$ is a homeomorphism. If $Z=[0,1]$, then each path in $C(X,Y)$ defines a continuous mapping $F: X \times [0,1] \to Y.$

\textbf{Definition.}   \index{homotopy}
Continuous mappings $f_0, f_1 \in C(X,Y)$ are called \textit{homotopic} $(f_0\sim f_1)$ if there exists such a continuous mapping $F: X \times [0,1]\to Y$ that $F (x,0)=f_0(x)$; $F(x,1)=f_1(x)$ for all $x\in X$. The continuous mapping $F$ is called a \textit{homotopy} between the mappings $f_0$ and $f_1$. The notation $f_t(x)=F(x,t)$ is also used.

\textbf{Examples.}
1) The mapping $F(x,t)=tx$ is a homotopy between the constant mapping to the origin and the identity mapping on $\mathbb{R}^n$.

2) If the topological space $Y$ is path-connected, then any two continuous mappings

3) Two mappings $f, g$ from an arbitrary space $X$ to an interval (a convex set) are homotopic to each other:

$$F(x, t) = (1-t) f(x)+ t g(x).$$

\begin{theorem} (Borsuk's Homotopy Extension Theorem) If $K$ is a CW complex, and $A$ is a subcomplex, then for any topological space $X$ and any continuous map $g: K \to X$ and homotopy $f_t: A \to X$ with $f_0=g | A$, there exists a homotopy $g_t: K \to X$, such that $f_t=g_t|A$ and $g_0=g$. \end{theorem}

\textbf{Corollary 1.} If $K$ is a CW complex, and $A$ is a subcomplex that can be contracted to a point, then $K/A$ is homotopically equivalent to $K$.

.

\textbf{Corollary 2.} If $(Y, A)$ is a pair of CW complexes and two attaching maps $f, g: A \to X$ are homotopic, then spaces $X\cup_fY$ and $X\cup_gY$ are homotopically equivalent.

\begin{theorem} (Cellular Approximation Theorem) Any continuous map $f: K \to L$ from one CW complex to another is homotopic to a cellular map $g$ ($g(e^k)\subset L^k$ for each cell $e^k$ of space $K$). \end{theorem}

\textbf{Problem.}     A CW complex is finite if and only if it is compact.  

\textbf{Problem.}     Prove that a CW space is connected if and only if its 1-skeleton is path-connected.  

\textbf{Problem.}     Prove that every connected CW complex is homotopically equivalent to a CW complex with a single vertex. (There exists a spanning tree where the vertices are all 0-cells, and the edges are some 1-cells).

The concept of homotopy defines an equivalence relation among the elements of the set $C(X,Y).$

\textbf{Definition.}  
The equivalence classes of homotopic mappings are called homotopy classes. The set of homotopy classes of continuous mappings from $X$ to $Y$ is denoted by the symbol $\pi(X,Y)$, that is, $$\pi(X,Y)=C(X,Y)/\sim.$$ For locally compact Hausdorff spaces, $\pi(X,Y)$ is the set of components of path-connectedness of $C(X,Y)$.

\textbf{Definition.}  
\index{retractable space}
A topological space is called retractable if the identity mapping $Id_X: X \to X$ is homotopic to a constant mapping (a mapping from $X$ to a point $x_0 \in X$). The homotopy between them is called the retraction of the space $X$ to the point $x_0$.
  
\textbf{Examples.} The spaces $I^n$, $D^n$ are retractable.

\textbf{Definition.}  \index{homotopy equivalence}
A mapping $f \in C(X,Y)$ is called a \textit{homotopy equivalence} of the spaces $X$ and $Y$ if there exists a mapping $g \in C(Y,X)$ such that $g \circ f \sim Id_X$; $f \circ g \sim Id_Y$. In this case, the mappings $f$ and $g$ are called \textit{homotopically inverse} (to each other) mappings. The class of homotopically equivalent spaces is called a \textit{homotopy type}. The notation for spaces is $X \sim Y$.

An example of homotopically equivalent spaces: a point and $\mathbb{R}^n$. It is evident that homeomorphic spaces are homotopically equivalent.

\textbf{Definition.}  
We say that a homotopy $F$ between continuous mappings $f, g: X \to Y$ is \textit{constant} on the set $A \subset X$ if for all $x \in A, t_1, t_2 \in [0, 1]: F(x, t_1) = F(x, t_2)$. That is, the images of points in the set $A$ remain fixed when the parameter $t$ changes. It is clear that the relation of homotopy that is constant on the set $A$ is also an equivalence relation.

\textbf{Problem.} 
Prove that the Möbius strip $Mo$ and the cylinder $S^1 \times [0,1]$ are not homeomorphic, but are homotopically equivalent.
 
\textbf{Problem.} 
Prove that the torus with a hole and the Klein bottle with a hole are not homeomorphic, but are homotopically equivalent.
 
\textbf{Problem.} 
Prove that $S^2 \setminus S^0$ is homotopically equivalent to $S^1$.
 
\textbf{Problem.} 
Prove that $S^3 \setminus S^1$ is homotopically equivalent to $S^1$.
 
\textbf{Problem.} 
Prove that $S^2 / S^0$ is homotopically equivalent to $S^1 \vee S^2$.
 
\textbf{Problem.} 
Prove that $S^2 / S^1$ is homotopically equivalent to $S^2 \vee S^2$.
 
\textbf{Problem.} 
Prove that every finite connected graph (1-complex) is homotopically equivalent to a graph (1-complex) with one vertex.
 
\textbf{Problem.} 
Prove that a space $X$ is contractible if and only if for any space $Y$, every continuous map $f: X \to Y$ is homotopic to a constant map if and only if every continuous map $f: Y \to X$ is homotopic to a constant.
 
\textbf{Problem.} 
Prove that the composition of homotopy equivalences is a homotopy equivalence. Show that a homotopy equivalence is an equivalence relation.
 
\textbf{Problem.} 
Prove that the infinite-dimensional sphere $S^{\infty}$ can be contracted to a point.

\section{Fundamental Group}

Recall that a path in a topological space $X$ is a continuous map $\alpha: I \to X$, where $I = [0,1]$. The point $x_0 = \alpha(0)$ is called the initial point, and $y_0 = \alpha(1)$ is called the terminal point of this path. The expressions are also used: $\alpha$ is a path from $x_0$ to $y_0$, the path $\alpha$ connects the points $x_0$ and $y_0$.

\textbf{Definition.}   \index{} A path whose initial and terminal points coincide is called a loop. The set of all loops in the space $X$ starting at the point $x_0$ is denoted by $\Omega(X, x_0)$.   

In the space $\Omega(X, x_0)$, we introduce the following notations:

1) \textit{Composition of loops} $u \circ v: I \to X$ starting at $x_0$ is defined by the condition
$$(u \circ v)(t) = \left\{
\begin{array}{c}
u(2t), \ \ \ \ \ \ \text{if} \ \ t \in [0, 1/2];\\
v(2t-1), \ \ \text{if} \ \ t \in [1/2, 1].
\end{array}
\right.$$
That is, it is a path that first traverses the loop $u$ and then the loop $v$.

2) A \textit{constant loop} at the point $x_0 \in X$ is defined as the constant map $u_0: I \to X;$ where $u_0(t) = x_0$ for all $t \in I$.

3) The \textit{inverse loop} to the loop $u^{-1}: I \to X$ is defined by the condition $u^{-1}(t) = u(1-t)$ for all $t \in I$, meaning it is a path that traverses the loop $u$ in the reverse direction.

\textbf{Definition.}   Two loops $u,v \in \Omega (X, x_0)$ are called \textit{equivalent} (denoted $u \sim v$) if there exists a homotopy between them that is constant at the endpoints 0 and 1 of the interval I.

\textbf{Definition.}   \index{fundamental group} The \textit{fundamental group} of the space $X$ relative to the point $x_0$ is the set of equivalence classes of loops $\Omega(X, x_0)/ \sim$, along with the operation of multiplication defined on it. This group is denoted $\pi_1(X, x_0)$.

\begin{theorem} 
The fundamental group of the circle $\pi_1(S^1, s_0)$ is isomorphic to the infinite cyclic group. It is generated by the equivalence class of the loop $u(t) = (\cos 2\pi t, \sin 2\pi t), \ 0 \le t \le 1.$
\end{theorem}

From this theorem, the following consequences can be derived:

1) The homeomorphic image of a circle (and any homotopic image to it) cannot be extended to a continuous mapping of a 2-disk.

2) (Brouwer's Fixed Point Theorem) Every continuous mapping of a disk onto itself has a fixed point.

3) Paths lying in a square and connecting opposite sides intersect.

4) (Jordan's Theorem) The complement of a closed continuous curve without self-intersections in the plane has two components of connectedness.

\textbf{Calculation of the Fundamental Group of a CW Complex.} Since the fundamental groups of homotopically equivalent spaces are isomorphic, without loss of generality, we can consider the CW complex $K$ to be connected with one vertex (0-cell) $x_0=e^0$.

1) The 1-skeleton $K^1$ is a bouquet of circles. Its fundamental group $\pi_1(K^1)$ is a free group $F(a_j)$ generated by $\{a_j\}$, which are defined by the 1-cells with chosen orientations on them.

2) Let’s compute $\pi_1(K^2, x_0)$. For each 2-cell $e_i^2$, the attachment map $f_i: \partial D^2 \to K^1$, together with a fixed point $y_0$ ($f_i(y_0)=x_0$) and the orientation of $D^2$, defines an element $R_i \in \pi_1(K^1, x_0)$ (a word composed of the letters $a_j^{\pm 1}$), which is trivial in $\pi_1(K^2, x_0)$. If $e^0 \cap f_i(\partial D^2) = \emptyset$, then $R_i=1$.
Thus, $\pi_1(K^2, x_0)$ has the representation $$(a_1, a_2,\ldots | R_1, R_2,\ldots ),$$ that is, the fundamental group $\pi_1(K^2, x_0)$ is generated by $a_j$ and the relations $R_i$ as the quotient group $F(a_j)/N(R_i)$.

3) The inclusion $K^2 \to K$ induces a map of the fundamental groups $\pi_1(K^2, x_0) \to \pi_1(K, x_0)$, which is an isomorphism.

\textbf{Example.} 
The fundamental group of a closed oriented surface of genus $g$ (given by the standard gluing of a polygon) is generated by $2g$ elements $a_i, b_i, i=1,\ldots , g$, which are related by a single relation $$a_1b_1a_1^{-1}b_1^{-1}a_2b_2a_2^{-1}b_2^{-1} \ldots a_gb_ga_g^{-1}b_g^{-1}.$$ In particular, the fundamental group of the two-dimensional torus is a free abelian group with two generators. When $g > 1$, the fundamental group of a closed oriented surface is non-abelian.

\textbf{Example.} 
The fundamental group of a closed non-orientable surface with $n$ Möbius strips is generated by $n$ elements $c_1,c_2, \ldots , c_n$, which satisfy the single condition $$c_1c_1c_2c_2\ldots c_nc_n=1.$$ In particular, the fundamental group of the projective plane is a cyclic group of order two. When $n > 1$, the fundamental group of a closed non-orientable surface is non-abelian.

\begin{theorem}
(Seifert–van Kampen) Let $X$ be a CW complex, $U$ and $V$ be CW subspaces, $X = U \cup V$, $W = U \cap V$, $x \in W$, $W$ is connected, $G = \pi_1(W, x)$. Then
$$\pi_1(X, x) = \pi_1(U, x) *_G \pi_1(V, x).$$
\end{theorem}
Here, the notation $\pi_1(U, x) *_G \pi_1(V, x)$ denotes the amalgamated product of the groups $\pi_1(U, x)$ and $\pi_1(V, x)$, which is a group whose generators are the generators of $\pi_1(U, x)$ and $\pi_1(V, x)$, and whose relations include the relations from $\pi_1(U, x)$ and $\pi_1(V, x)$, as well as the relation $$i_*(a_k) = j_*(a_k),$$ where $i:W \to U,$ $j:W \to V$ are inclusion maps, and $\{a_k\}$ is a system of generators of the group $G$. It is assumed that all groups have a finite number of generators and relations.

The theorem also holds in the case of any connected space $X$, provided that $U$, $V$, and $W$ are connected open subsets.

\textbf{Problem.} 
Find the fundamental groups of the following spaces:

1) torus; 2) torus with a hole; 3) torus with two holes; 4) sphere with three holes; 5) Klein bottle; 6) Möbius strip with a hole; 7) Klein bottle with two holes.

\textbf{Problem.} 
Prove that a loop is homotopic to 0 if and only if the mapping of the circle that defines it extends to a mapping of the 2-dimensional disk.

\textbf{Problem.} Prove that any continuous mapping of a circle onto its circle is homotopic to a constant mapping.

\section{Literature Review}















Many works have been dedicated to the study of topological properties of functions and dynamical systems. Among them, we  highlight the works of the Kyiv topology school: 
scientific articles by the author on function topology \cite{Prishlyak1993, prishlyak1998, prishlyak1999equivalence, prishlyak2000conjugacy, prishlyak2001conjugacy, prishlyak2002topological1, prishlyak2002morse, prishlyak2003regular}, vector fields \cite{Prish1997vec, Prishlyak2001, prishlyak2002morse1, Prishlyak2002, prishlyak2003topological, prishlyak2003sum, prishlyak2005complete, Prishlyak2007}, and other geometric objects \cite{Prishlyak1994, prishlyak1997graphs, Prishlyak1999}, as well as the works of his students: K. Myshchenko \cite{prishlyak2007classification}; N. Budnytska \cite{Bud2008knu}; D. Lychak \cite{lychak2009morse}; A. Bondarenko \cite{Bond2012mfat}; O. Vyatychaninova \cite{VyatP2013Mol}; Bohdana Hladysh \cite{Hladysh2016, hladysh2017topology,
Hladysh2019}; A. Prus \cite{Prishlyak2017, prishlyak2020three, Prishlyak2021}; V. Kiosak \cite{KPL2022}; S. Bilun \cite{bilun2002closed}
V. Lisikevich \cite{LisP2013KNU}, I. Ivaniuk \cite{IvanPrish2014-func-deform3, IP2015knu}, D. Skotchko \cite{
PS2016-F-atoms}, M. Loseva \cite{Losieva2017, Prishlyak2019, Prishlyak2020}, Z. Kybalko \cite{Kybalko2018}, K. Khatamian \cite{Hatamian2020}.






\end{document}